\newfont{\blackboard}{msbm10}
\def\bbf#1{\mbox{\blackboard #1}}
\def\goth#1{\mathfrak #1}
\def\Fd{{\mathbb D}}

\documentclass{article}
\usepackage{amsfonts}
\begin{document}

\newtheorem{theorem}{Theorem}
\title{The Differential and Functional Equations for a Lie Group
Homomorphism are Equivalent}
\author{George Svetlichny}
\maketitle

\begin{abstract}
I prove the ``folklore" result that the functional equation for a Lie
group homomorphism can be solved by solving the corresponding
differential equation.
\end{abstract}

\section{Introduction}
      The simplest example of our result concerns the functional equation:
\[f(x+y) = f(x)f(y)\]
As is well known, the measurable solutions with \(f(0) = 1\) are
exponentials: \(f(x) = e^{kx}\). If we differentiate the equation in
\(y\) and
set \(y = 0\) we get:
\[f'(x) = f'(0)f(x)\]

The solutions of this {\em differential\/}  equation are the same exponential
functions with \(k = f'(0)\). A more complicated example occurred
\cite{twodimref}  when
the I was faced with the problem of finding all  lower triangular
\(2\times 2\) matrices \(M(X)\) depending on a vector \(X\)   satisfying
\(M(0) = I\) and the functional equation:
\[ M(X+M(X)Y) = M(X)M(Y)\]
The (differentiable) solutions were found to coincide with the
solutions of the system of differential equations resulting from
differentiating with respect to \(Y\) and setting \(Y = 0\). Now it is far
from common that a functional equation imply a differential equation,
and even less common that all solutions of the differential equation
be solutions of the original functional equation. What makes the trick
work in the above cases is the fact that both equations can be
construed as stating a Lie group homomorphism, and for these, as we
shall shortly show, the functional and differential equations are
equivalent, at least in a neighborhood of the identity, provided we are interested in differentiable solutions.\footnote{ For the functional equation it's probably enough to just require
measurability.} I have seen this result used on many occasions but have never seen a proof and so decided do provide one myself.
\section{The Result}
 All manifolds and maps are considered
to be \({\cal C}^\infty\).
   Let \(G\) be a Lie group, \(U\) a neighborhood of the identity \(e\in
   G\), and
\(M\) a manifold of dimension \(m\). Let \(p:U \to M\) be a fibration and
 let \(x_0
 = p(e)\). Consider now the problem of finding in a neighborhood of
 \(x_0\) a
 section  \(\sigma\) of \(p\) such that:
\begin{eqnarray}\label{eq:prohomid}
           \sigma(x_0) & = & e  \\  \label{eq:prohomprod}
 \sigma(p(\sigma(x) \sigma(y))) & = &  \sigma(x) \sigma(y)
\end{eqnarray}

    Our two examples are instances of this. For the first example  let
\({\bbf R}^*\) be the multiplicative group of positive reals and take
\(G\) to  be the
direct product  \({\bbf R}^* \times {\bbf R}\)  with \(p\)
the projection onto the second
factor; \(f\) is then the first component of  \(\sigma\).
For the second example
let \(T\)  be the group of invertible lower triangular \(2 \times 2\)
 matrices
and take \(G\) to be \(T \times {\bbf R}^2\) with the product law:
\[(M,X)\cdot(N,Y) = (MN,X+MY),\]
and \(p\) be the projection on the second factor; as before the desired
function is the first component of  \(\sigma\).

    Now what (\ref{eq:prohomid}-\ref{eq:prohomprod}) effectively say is that a neighborhood
    of \(x_0\)  is to
acquire a local Lie group structure with the product given by:
\[x,y \mapsto  p( \sigma(x) \sigma(y))\]
and that  \(\sigma\) is to provide a homomorphism of this structure into
\(G\).
Thus the image of  \(\sigma\) is a local Lie subgroup of \(G\) of
dimension \(m\)
and transversal to the fiber \(p^{-1}(x_0)\). Reciprocally, any local Lie
subgroup with these two properties is a solution by taking  \(\sigma\) to be
the inverse of \(p\) restricted to this subgroup. The problem can be
further reduced to an algebraic one. The tangent space \({\goth g}\)   to \(G\) at
\(e\)
(the Lie algebra of \(G\)) has a natural vertical subspace (not
necessarily a subalgebra) provide by the tangent space to the fiber.
The germs of \(m\)-dimensional local Lie subgroups transversal to the
fiber are now in one to one correspondence to the \(m\) dimensional Lie
subalgebras of \({\goth g}\)  that, as subspaces, are transversal to
the vertical
subspace. This in principle resolves  the existence and uniqueness
problem for germs at \(x_0\)  of solutions to (\ref{eq:prohomid}-\ref{eq:prohomprod}). One hasn't though
resolved the {\em practical\/}  problem of finding such solutions explicitly,
for even if the algebraic problem is solved one has to exponentiate
the Lie subalgebra to find the map and this may not be a trivial task.
Fortunately, the algebraic and exponentiation problems can be avoided
by appealing to the differential equation that results from
(\ref{eq:prohomprod}). For
ease of notation let \(\mu\) denote the group product of \(G\). If we
differentiate  (\ref{eq:prohomprod}) with respect to \(y\) and set \(y =
x_0\)  we obtain the
following differential equation for  \(\sigma\):

\begin{eqnarray}\label{eq:homdif}
 \Fd \sigma(x) \cdot \Fd p( \sigma(x)) \cdot \Fd_2 \mu( \sigma(x),e) \cdot \Fd
\sigma(xo) &=& \Fd_2 \mu( \sigma(x),e) \cdot \Fd \sigma(x_0)
\end{eqnarray}
where \(\Fd\) denotes the Frechet derivative. Now obviously any solution to
(\ref{eq:prohomid}-\ref{eq:prohomprod}) necessarily satisfies (\ref{eq:homdif}), what is remarkable is that the converse is true.
\begin{theorem}
      Any local solution of (\ref{eq:homdif}) such that  \(\sigma(x_0) = e\)
       is a (local) solution of (\ref{eq:prohomid}-\ref{eq:prohomprod}).
\end{theorem}
     {\em Proof:}  By our previous discussion it's enough to show that the
image \(S\) of  \(\sigma\) is a local Lie subgroup of \(G\). Apply both sides of (\ref{eq:homdif})
to a tangent vector  \(\eta\) at \(x_0\) . The left hand side is of the form
\(\Fd\sigma(x) \cdot \eta\) and is thus a tangent vector to \(S\) at
\(\sigma(x)\). The right-hand side is
\(\Fd_2 \mu( \sigma(x),e) \cdot \Fd \sigma(x_0) \cdot \xi\). Now \(\Fd
\sigma(x_0) \cdot \xi\) is a tangent vector to \(S\)
at \(e\), and by picking  \(\xi\) appropriately any such tangent
vector can  be
so given. On the other hand \(\Fd_2 \mu( \sigma(x),e)\) is the tangent map of left
multiplication by  \(\sigma(x)\). Thus a consequence
of (\ref{eq:homdif}) is that the left
translate by  \(\sigma(x)\) of a tangent vector to \(S\) at \(e\)
is tangent to \(S\) at
\(\sigma(x)\), in other words: those left-invariant vector fields that are
tangent to \(S\) at \(e\) are tangent to \(S\) (at all other points).
Because of
this tangency to the same submanifold, these vector fields are in
involution and so form a Lie subalgebra. Associated to this subalgebra
is a unique germ of a local Lie subgroup. Since any such subgroup has
a neighborhood of the identity covered by exponentiations of the
tangent left-invariant vector fields, we see that the subgroup germ
coincides with the germ of \(S\) at \(e\).  Q.E.D

\section{Acknowledgments}

This research received partial financial support from the Conselho Nacional de Desenvolvimento Cient\'{\i}fico e Tecnol\'ogico (CNPq) and the Funda\c{c}\~ao de Amparo \`a Pesquisa do Estado do Rio de Janeiro  (FAPERJ).

\end{document}